\documentclass[a4paper,10pt]{article}
\usepackage{amsmath,amsthm,amssymb,amsfonts,epsfig,epstopdf,titling,url,array,enumerate,cite,mathrsfs,upgreek,authblk,scrextend, graphicx, float}
\usepackage[bindingoffset=0.2in,left=1in,right=1in,top=1in,bottom=1in,footskip=.25in]{geometry}
\theoremstyle{plain}
\newtheorem{thm}{Theorem}[section]
\providecommand{\keywords}[1]{\begin{addmargin}[28pt]{28pt}\noindent\textbf{Keywords:} #1 \end{addmargin}}

\newtheorem{ppn}[thm]{Proposition}
\theoremstyle{definition}
\newtheorem{dfn}[thm]{Definition}
\newtheorem{eg}[thm]{Example}
\newtheorem{rem}[thm]{\textit{Remark}}
\providecommand{\ams}[1]{\begin{addmargin}[28pt]{28pt}\noindent\textbf{Mathematics Subject Classification:} #1\end{addmargin}}

\title{\textbf{Some fixed point theorems in generalized parametric metric spaces and applications  to ordinary differential equations}}
\author[1]{Abhishikta Das}
\author[2]{Hijaz Ahmad}
\author[3]{T. Bag\footnote{corresponding author}}
\affil[1,3]{Department of Mathematics, Siksha-Bhavana,		Visva-Bharati, Santiniketan-731235, West-Bengal, India }
\affil[2]{Section of Mathematics, International Telematic University Uninettuno, Corso Vittorio Emanuele II, 39,00186 Roma, Italy 
	\authorcr 	E-mail\textsuperscript{1}: abhishikta.math@gmail.com 
	\authorcr   E-mail\textsuperscript{2}: hijaz555@gmail.com 
	\authorcr   E-mail\textsuperscript{3}: tarapadavb@gmail.com}
\affil[1]{Orcid Id: 0000-0002-2860-424X }
\affil[2]{Orcid Id: 0000-0002-5438-5407}
\affil[3]{Orcid Id: 0000-0002-8834-7097 }

\date{}

\begin{document}
	\maketitle
	
	\begin{abstract}
		\noindent
		The objective of this work is the construction of `Boyd-Wong fixed point theorem' in the setting of  generalized parametric metric space and discussion its application on existence criteria of solutions to a second order initial value problem. Also an analogue of `Banach type fixed point theorem' of generalized parametric metric space is proved and its application on the existence of a solution of first order periodic boundary value differential equation is examined.
		
	\end{abstract}
	
	\keywords{Generalized parametric metric space, contraction mapping,  Boyd-Wong fixed point theorem. } 
	\ams{47H10, 54H25, 34B15.}

	\section{Introduction and Preliminaries}
	The study of  fixed point theory in metric space is one of the most important, interesting and applicable branch among all the branches of mathematical theories. As in 1906, mathematician M. Frechet \cite{2}  proposed the concept of metric space, several mathematicians used to involve themselves with the study of metric space and metric fixed point theory.  As a result we became enriched by the notion of several types of generalized distance functions(viz. b-metric\cite{3}, S-metric\cite{4}, cone-metric\cite{5}, $\phi$-metric\cite{6} etc.). Most of these generalizations  were  done by changing the `triangle inequality' of metric axioms or some times by introducing more co-ordinates in the domain of the distance function.      \\
	Recently, in 2014, N. Hussian et al.\cite{7} introduced  parametric metric space, a new generalized metric space,  where they used a non-negative parameter `t' in the definition and studied fixed point results on this new setting. After that some researchers attempted to  generalize the thought of parametric metric in different approaches(please see \cite{12,13,14,15,16} etc.). Analyzing such hybrid definitions of parametric metric, Das and Bag\cite{1} recently introduced `generalized parametric metric space'. \\
	In the subsequent results we recall the notion  of generalized parametric metric space  and some  properties which are required for this article. 
	\begin{dfn} \cite{1} \label{dfn1}
		Let $ A = \mathbb{R}_{\geq 0} $ and  define a binary operation $ o : A \times A \rightarrow \mathbb{R}_{\geq 0}  $  satisfying the following conditions: 	  
		\begin{enumerate}[(i)]
			\item $ \alpha ~ o ~ 0 = \alpha ~~ $ 
			\item $  \alpha \leq \xi $ implies $ \alpha ~ o ~ \zeta \leq \xi ~ o ~ \zeta ~~ $ 
			\item $ \alpha ~ o ~ \xi =  \xi ~ o ~ \alpha ~~ $ 
			\item $ \alpha ~ o ~ ( \xi ~ o ~ \zeta ) = (  \alpha ~ o  ~ \xi ) ~ o ~ \zeta ~~ $ 
		\end{enumerate}	
		for all $ \alpha, ~ \xi, ~ \zeta \in A  $. 
	\end{dfn}

	Some examples of such $ `o $' are: 
	\begin{enumerate}[(i)]
		\item $ \xi ~ o ~ \zeta = \max \{\xi ,  \zeta \} $
		\item $\xi ~ o ~ \zeta = \xi ~ + ~ \zeta  $
	\end{enumerate}

	\begin{dfn} \cite{1} \label{dfn2}
		Let $ \{ a_n \} $ and $  \{ b_n \} $ be two sequences in $ \mathbb{R}_{\geq 0} $ which converges to $ a $ and $b $ respectively. If $  \{a_n  o b_n \} $ converges to $ a o  b $ then 
		$ ` o $' is said to be continuous.
	\end{dfn}

	In the following we recall   some other  axioms of $ ` o $'. \\
	(v)  $ ` o   $' is a continuous function. \\
	(vi) $  \alpha < \xi $ and $ \zeta < \mu $  implies $ \alpha ~ o ~ \zeta < \xi ~ o ~ \mu, ~ $   for all $ \alpha, ~ \xi, ~	\zeta, ~ \mu \in \mathbb{R}_{\geq 0}. $ \\
	(vii) $ \xi ~ o ~ \xi \geq \xi, ~ $ for all $ \xi \in \mathbb{R}_{\geq 0}.~~$

	\begin{dfn}\cite{1} \label{dfn3}
		A generalized parametric metric space $ ( \Im, \rho, o ) $ is a triple where $ \Im \neq \phi $,  \\
		$ \rho : \Im \times \Im \times (0, + \infty) \rightarrow [0, + \infty) $ is a function, called  generalized parametric metric  satisfying the following properties for all  $  x, \zeta, z \in \Im $ 
		\begin{enumerate}[($\rho1$)]
			\item $ ( \rho ( x, \zeta, r ) = 0 ~ \text{for all} ~ r > 0) ~ $  if and only if $ x = \zeta $;
			\item $  \rho ( x, \zeta, r ) = \rho ( \zeta, x, r ) ~~$ for all $ r > 0 $;
			\item  $  \rho (x, \zeta, r_1 + r_2 ) \leq  \rho (x, z, r_1) ~ o ~ \rho (\zeta, z, r_2 ) $ for all $ ~ r_1, r_2 > 0 $.
		\end{enumerate}	
	\end{dfn}

	\begin{eg} \cite{1} \label{eg1}
		Consider the function
		$$ \rho ( x, \zeta, r )= \frac{| x - \zeta |^p}{r} ~ \text{ for all} ~ x, \zeta \in \mathbb{R}, ~  r > 0 ~ \text{and}~ 0 < p < 1. $$ 
		Then $ (\mathbb{R}, \rho, + ) $ is a generalized parametric metric space.  
	\end{eg}

	\begin{ppn} \cite{1} \label{lma1}
		Every generalized parametric  metric $ \rho ( a,b, s ) $    on $ \Im $  is non-increasing function with respect to $ s $.
		
	\end{ppn}
	
	\begin{dfn} \cite{1} 
		Let   $ \{ x_n \}  $ be a sequence in a generalized parametric metric space $ ( \Im, \rho, o ) $. Then   $ \{ x_n \}  $ is said to be a 
		\begin{enumerate}[(i)]
			\item   convergent sequence if  there exists $ \zeta \in \Im $ such that for all $ l > 0 $,
			$ \rho ( x_n, \zeta, l ) \rightarrow  0 $ as $ n \rightarrow + \infty  $. 
			\item  Cauchy sequence if for all $ l > 0 $, $  \rho ( x_n, x_m, l ) \rightarrow  0 $ as $ m, n \rightarrow + \infty $.    
		\end{enumerate} 
	\end{dfn}

	\begin{ppn} \cite{1} 
		In a generalized parametric metric space $ ( \Im, \rho, o ) $ if  $ ` o  $' is continuous, then
		\begin{enumerate}[(i)]                       
			\item limit of  convergent sequence is unique;
			\item every convergent sequence in $ \Im $ is a Cauchy sequence.
		\end{enumerate}
	\end{ppn}

	\begin{dfn} \cite{1} 
		Let $ ( \Im, \rho , o ) $ be a	 generalized parametric metric space.   If each  Cauchy sequence in $ \Im $ is convergent and converges to some member in $  \Im $ then  $( \Im, \rho , o ) $ is called complete.
	\end{dfn} 
	
	Next theorem is the presentation of  `Banach contraction theorem'  in  generalized parametric metric space setting.

	\begin{thm} \cite{1}
		Let $ \mathcal{F} $ be a self mapping defined on a complete generalized parametric metric space $ ( \Im, \rho , o ) $  and $ ` o  $' be continuous. If for all $ x, \xi \in \Im $ and for all $ \alpha > 0 $, $ \mathcal{F} $ satisfies the following  condition 
		$$ \rho ( \mathcal{F} x, \mathcal{F} \xi, \alpha ) \leq c \rho (x, \xi, \alpha ) $$
		where $ 0 < c < 1 $, then $ \mathcal{F} $ has a unique fixed point in $ \Im $.
	\end{thm}

	We are conscious with the fact that, some of the metric fixed point results in generalized metric spaces are very obvious due to the metrizability of that space. But  proof of some non-linear type contractions, for example, the famous `Boyd-Wong fixed point theorem'\cite{8} is very interesting	for different generalized metric spaces and not necessarily followed by the result of classical metric space (for reference see \cite{new ref}). That's why we exercise on the construction of    Boyd-Wong type result in  generalized parametric metric space. We also apply the Boyd-Wong type fixed point result of generalized parametric metric space to  find a solution of the following second order initial value problem 
	\begin{equation} \label{B-W 1}
		\dfrac{d^2 \mathcal{Y}}{d x^2 } + w^2 \mathcal{Y} = \mathcal{G} ( x, \mathcal{Y}( x ) ),  ~~~ \mathcal{Y} ( 0 ) = l_1, ~ \mathcal{Y}' ( 0 ) = l_2
	\end{equation}
	where $ \mathcal{Y} \in C [ 0, S ] $, $ w ( \neq 0 ), ~ \mathcal{G} : [ 0, S ] \times \mathbb{R}^+ \rightarrow \mathbb{R} $  is a function and $  l_1, l_2 \in \mathbb{R} $. \\
	
	In another Section, we study the theory of fixed point  in partially ordered set.  In 2003,  Ran and  Reurings\cite{9}  proved an analogue result of Banach fixed point theorem\cite{10} in partially ordered sets and established it's various applications.  Being inspired by those developments, we established an analogue result of Banach type fixed point theorem in generalized parametric metric space\cite{1} to  partially ordered set. The key principal in this type of   theorem is that the contraction condition for the nonlinear monotone self mapping is supposed to hold only for the  elements that are comparable in the partial order relation of the  set.  Main focus of  this result is that under such conditions the conclusions of Banach type fixed point theorem of generalized parametric metric space still hold over partially ordered sets. \\
	
	Next we recall the following  definition related to  differential equation which is required in our main result. \\

	We consider periodic boundary value problem 
	\begin{equation} \label{periodic ode1}
		u'(  \mathcal{Y}) = \mathcal{F} ( \mathcal{Y}, u( \mathcal{Y} ) ), ~ \mathcal{Y} \in [0, S ], ~~ u ( 0 ) = u ( S ) 
	\end{equation}
	where $ S > 0 $ and $ \mathcal{F} :  [0, S ] \times \mathbb{R} \rightarrow \mathbb{R} $ is a continuous mapping.
	\begin{dfn} \cite{11}
		A solution for the system (\ref{periodic ode1}) is a function $ \alpha \in C [0, S ] $  satisfying the conditions of (\ref{periodic ode1}). \\
		A lower solution of (\ref{periodic ode1}) is a function $ \alpha \in C [0, S ] $ such that  $ 	\alpha'( \mathcal{Y} ) \leq \mathcal{F} ( \mathcal{Y}, \alpha ( \mathcal{Y} ) ), ~ \mathcal{Y} \in [0, S ], ~~ \alpha ( 0 ) \leq \alpha ( S ) $. \\
		An upper solution of (\ref{periodic ode1}) is a function $ \alpha \in C [0, S ] $ satisfying the reverse relation.
	\end{dfn}
	
	The presentation of this manuscript is as follows.\\
	In Section 2, Boyd-Wong type fixed point theorem and its application on finding  solution of second order initial value differential equation is established. Section 3 consists of Banach type fixed point theorem  in generalized parametric metric space over partially ordered set and its application to the existence criteria of solution for first order periodic boundary value differential equation.

	\section{Boyd-Wong type contraction in generalized parametric metric space and its application to differential equation }
	
	In this Section first we introduce the notion of Boyd-Wong type contraction in  generalized parametric metric space setting and then establish fixed point theorem for such type of mappings.  An application of this theorem for   existence and uniqueness criteria of solutions for  second order initial value differential equation is given. 
	
	\subsection{Boyd-Wong type contraction}
	
	\begin{dfn}
		A self mapping  $ \mathcal{T} $ on a 	generalized parametric metric space $ ( \Im, \rho, o ) $ is said to satisfy Boyd-Wong type contraction condition if for all $ \alpha, \zeta \in \Im $, $ \mathcal{T} $ satisfies 
		\begin{equation}\label{equ 1}
			\rho ( \mathcal{T} \alpha, \mathcal{T} \zeta, l ) \leq \phi ( \rho ( \alpha, \zeta, l ) ) ~ ~ \text{for all} ~l > 0 
		\end{equation}
		where $ \phi $ is  upper semi-continuous from the right, defined from  $ \mathbb{R}^+ $ to $ \mathbb{R}^+ $ and  satisfies $ \phi ( s ) < s ~ $ for all $ s > 0 $.
	\end{dfn}

	\begin{thm} \label{main theorem}
		Let $ ( \Im, \rho, o ) $ be a complete generalized parametric metric space and $ \mathcal{T} $ be a function on $ \Im  $ which satisfies  Boyd-Wong type contraction condition. 
		Then $ \mathcal{T} $ has a unique fixed point $ \alpha^* \in \Im $  and $ \{ \mathcal{T}^n \alpha \} $ converges to $ \alpha^* $ for each $ \alpha \in \Im $.
		\begin{proof}
			Let us choose $ x_0 \in \Im $ and  define $ \zeta_n( s ) = 
			\rho ( \mathcal{T}^{n-1} x_0, \mathcal{T}^n x_0, s ) $, for $ s > 0 $. \\
			Then  for all $ q \in \mathbb{N} $ and for each $ s > 0 $, we have 
			$$	\zeta_{q+1} ( s ) =  \rho ( \mathcal{T}^{q} x_0, \mathcal{T}^{q+1} x_0, s )  \leq \phi ( \rho ( \mathcal{T}^{q-1} x_0, \mathcal{T}^q x_0, s ) ) 	 = \phi ( \zeta_q ( s ) ) < \zeta_q ( s ). $$
			This shows that $ \{ \zeta_n ( s ) \} $ is a decreasing sequence for each $ s > 0 $ and bounded below. \\ Hence $ \{ \zeta_n ( s ) \} $ has a limit, say $ a ( s ) $ for each $ s > 0 $. \\
			Next, we claim that $ a ( s )= 0 ~ $ for all $ s > 0 $. \\
			If possible suppose,  there exists $ s_0 > 0 $ such that $ a ( s_0 ) > 0 $. \\
			Then $ \zeta_{n+1} ( s ) \leq \phi ( \zeta_n ) ( s ) $ and upper semi-continuity from the right of $ \phi $ gives, 
			$$ a ( s_0 ) \leq \underset{ \zeta_n ( s_0 ) \rightarrow a ( s_0 ) ^+}{lim} \phi ( \zeta_n ( s_0 ) ) \leq \phi ( a ( s_0 ) ). $$
			But this contradicts that $ \phi ( s ) < s ~ $  for all $ s > 0 $. \\
			Thus we get,
			$$  \underset{n \rightarrow + \infty}{lim} a_n ( s ) = \underset{n \rightarrow + \infty}{lim} \rho ( \mathcal{T}^{n-1} x_0, \mathcal{T}^n x_0, s ) = 0 ~  \text{for all} ~ s > 0. $$ 
			Next claim is that $ \{ \mathcal{T}^n x_0 \} $ is a Cauchy sequence in $ \Im $. \\
			For each  $ t > 0, ~ n \in \mathbb{N} $ and $ q = 1, 2, 3, \cdots $, we have
			\begin{align*}
				& 	\rho (  \mathcal{T}^{n+q} x_0, \mathcal{T}^n x_0, s ) \\
				& \leq \rho (  \mathcal{T}^{n+q} x_0, \mathcal{T}^{n+1} x_0, \frac{s}{2} ) ~ o ~ \rho (  \mathcal{T}^{n+1} x_0, \mathcal{T}^{n} x_0, \frac{s}{2} ) \\
				& \leq \rho (  \mathcal{T}^{n+q} x_0, \mathcal{T}^{n+2} x_0, \frac{s}{4} ) ~ o ~ \rho (  \mathcal{T}^{n+2} x_0, \mathcal{T}^{n+1} x_0, \frac{s}{4} ) ~ o ~ \rho (  \mathcal{T}^{n+1} x_0, \mathcal{T}^{n} x_0, \frac{s}{2} ) \\
				& \cdots \hskip 10 pt \cdots \\
				& \leq \rho (  \mathcal{T}^{n+q} x_0, \mathcal{T}^{n+q-1} x_0, \frac{s}{2^{q-1}} ) ~ o ~  \rho (  \mathcal{T}^{n+q-1} x_0, \mathcal{T}^{n+ q -2} x_0, \frac{s}{2^{q-1}} ) ~ o ~ \cdots ~ o  ~ \rho (  \mathcal{T}^{n+1} x_0, \mathcal{T}^{n} x_0, \frac{s}{2} ) 
			\end{align*} 
			This yields, 	
			$$ \underset{n \rightarrow + \infty}{lim} \rho ( \mathcal{T}^{n+q} x_0, \mathcal{T}^n x_0, s  ) = 0 ~~ \text{for all} ~ s > 0 ~ \text{and} ~ q = 1, 2, \cdots. $$
			Hence $ \{ \mathcal{T}^n x_0 \} $ is a Cauchy sequence in $ \Im $ and thus the sequence converges to some element, say $ u^* $ in $ \Im $. \\
			Again for all $ ~ s > 0 $, we have
			$$ \rho ( u^*, \mathcal{T} u^*, s )   = \underset{n \rightarrow + \infty}{lim} \rho ( \mathcal{T}^n x_0, \mathcal{T} u^*, s )  \leq 	\underset{n \rightarrow + \infty}{lim} \phi ( \rho ( \mathcal{T}^{n-1} x_0,  u^*, s ) ) < \underset{n \rightarrow + \infty}{lim} \rho ( \mathcal{T}^{n-1} x_0,  u^*, s ) = 0. $$ 
			This implies $	\rho ( u^*, \mathcal{T} u^*, s )   = 0 ~~ \text{for all} ~ s > 0 ~~ $ i.e. $ ~  \mathcal{T} u^* = u^* $. \\
			Moreover,   $ u^* $ is the only fixed point for $ \mathcal{T} $. \\
			Since, if  there exists another fixed point $  u^* $   of $ \mathcal{T} $, then
			$$	\rho ( u, u^*, s ) = \rho ( \mathcal{T} u, \mathcal{T} u^*, s ) \leq \phi ( \rho ( u, u^*, s ) ) < \rho ( u, u^*, s ) ~~ \text{for all}  ~ s > 0 ~ \text{implies} ~ u = u^*. $$
			This completes the proof.
		\end{proof}
	\end{thm}
	
	\begin{eg} \label{eg 2}
		Let $ \Im = [-1, 1] $. We consider the generalized parametric metric space $ ( \Im, \rho, o ) $ where $ \rho $ is defined by $ \rho ( a, \xi, l ) = \dfrac{\sqrt{ | a - \xi | } }{ l } $ for all $ l > 0 $ and $ a, \xi \in \Im $ (for details please see \cite{1}). \\
		Let $ \phi ( \mu ) = \frac{ \mu }{ 2 } ~ $  for all $ \mu \in  \mathbb{R}^+  $. \\
		Then clearly $ \phi $ is upper semi-continuous from the right and $ \phi ( \mu ) < \mu $ for all  $ \mu \in  \mathbb{R}^+  $. \\
		Next we define a function $ {T} $ on $ \Im $ by $ {T} ( \mu ) = \frac{\mu}{16} ~ $ for all $ \mu \in \Im $. \\
		Then,
		$$ \rho ( {T}  a, {T}  \mu, s ) = \dfrac{\sqrt{| \frac{a}{16} - \frac{ \mu }{16}|}}{s} = \frac{1}{4} \dfrac{\sqrt{ | a - \mu | }}{s} ~~ \text{and} ~ ~ \phi( \rho ( a, \mu, s ) ) = \phi (\dfrac{\sqrt{ | a  - \mu | }}{ s }) = \frac{1}{2}\dfrac{\sqrt{ | a - \mu | } }{s} $$ 
		for all $ s >  0 $  and   $ a, \mu \in X $. \\
		Therefore,
		$$ \rho ( {T}  a, {T}  \mu, t ) <  \phi( \rho ( a, \mu, t ) ) ~  \text{for all} ~ t > 0 ~ \text{and} ~ a, \mu \in \Im. $$
		Hence by the Theorem \ref{main theorem}, $ {T}  $ is a mapping having unique fixed point 
		and here which is  $ 0 $. 

		\section*{Illustration of the Example  \ref{eg 2} by graphically and numerically:} 
		Figure 1 shows the illustration of  Example \ref{eg 2} on 2D view, in which the variation of $ \rho ( {T} a, {T} \mu, t ) $ and $ \phi( \rho ( a, \mu, t ) ) $   as a function of $ x $ where $ x = | a - \mu |, $ for all $ a, \mu \in \Im $ with fixed values of $ t = 1 $ is shown as a red and blue  colored curves  and denoted as $ H (x,1) $ and $ G ( x, 1 ) $ respectively. \\
		On the other hand, Figure 2 shows the illustration of Example \ref{eg 2} on 2D view, in which the variation of $ \rho ( {T} a, {T} \mu, t ) $ and $ \phi( \rho ( a, \mu, t ) ) $   as a function of $ t $ and fixed values of $ x = 1 $ 	where $ x = | a  - \mu |, $  for all  $ a, \mu \in \Im $ is shown as a red and blue  colored curves  and denoted as $ H ( 1, t ) $ and $ G ( 1, t ) $ respectively. \\
		\begin{figure}
			[h!]
			\centering
			(a) \includegraphics[width=6cm,height=5.5cm]{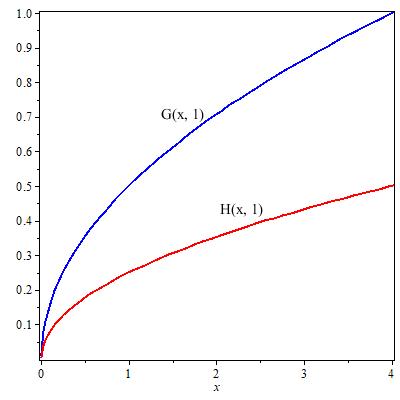}
			%
			(b) \includegraphics[width=6cm,height=5.5cm]{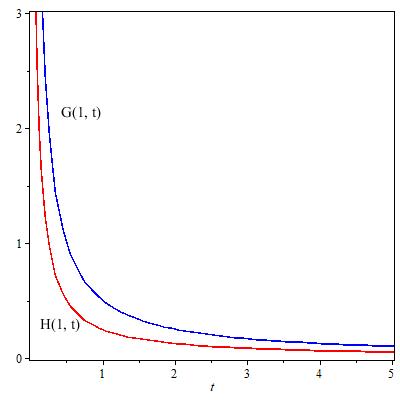}
			\caption{(a) $ H (x,1) $ vs $ G (x,1) $ $ \hskip 5 pt $ (b) $ H ( 1, t ) $ vs $ G ( 1, t ) $}
		\end{figure}

		The following tables(Table 1 and 2) numerically justifies the inequality (\ref{equ 1}) for the Example \ref{eg 2} i.e shows the variation of $ \rho ( {T} a, {T} \mu, t ) $ and $ \phi ( \rho ( a, \mu, t ) ) $  which is observed for both the curves with respect to  the values of $  x $ where $ x = | a - \mu |, ~  a, \mu \in \Im $ and for different fixed ranges of $ t $ and vice versa. \\
		Table 1 shows the variation between      $ \rho ( {T} a, {T} \mu, t ) $ and $ \phi ( \rho ( a, \mu, t ) ) $ as a function of $ x $   where $ x = | a - \mu |, ~  a, \mu \in \Im $ with relative to two different ranges of $ t $ which are taken as $ t = 1 $ and $ t = 30 $. \\
		Table 2 shows the variation between      $ \rho ( {T} a, {T} \mu, t ) $ and $ \phi ( \rho ( a, \mu, t ) ) $ as a function of $ t $ with relative to two different ranges of $ | a - \mu | \in [ 0, 2 ] $ which are taken as $ 0.5 $ and $ 1.8 $.\\

		\begin{table}[h!]
			\centering
			\includegraphics[width=10cm,height=5cm]{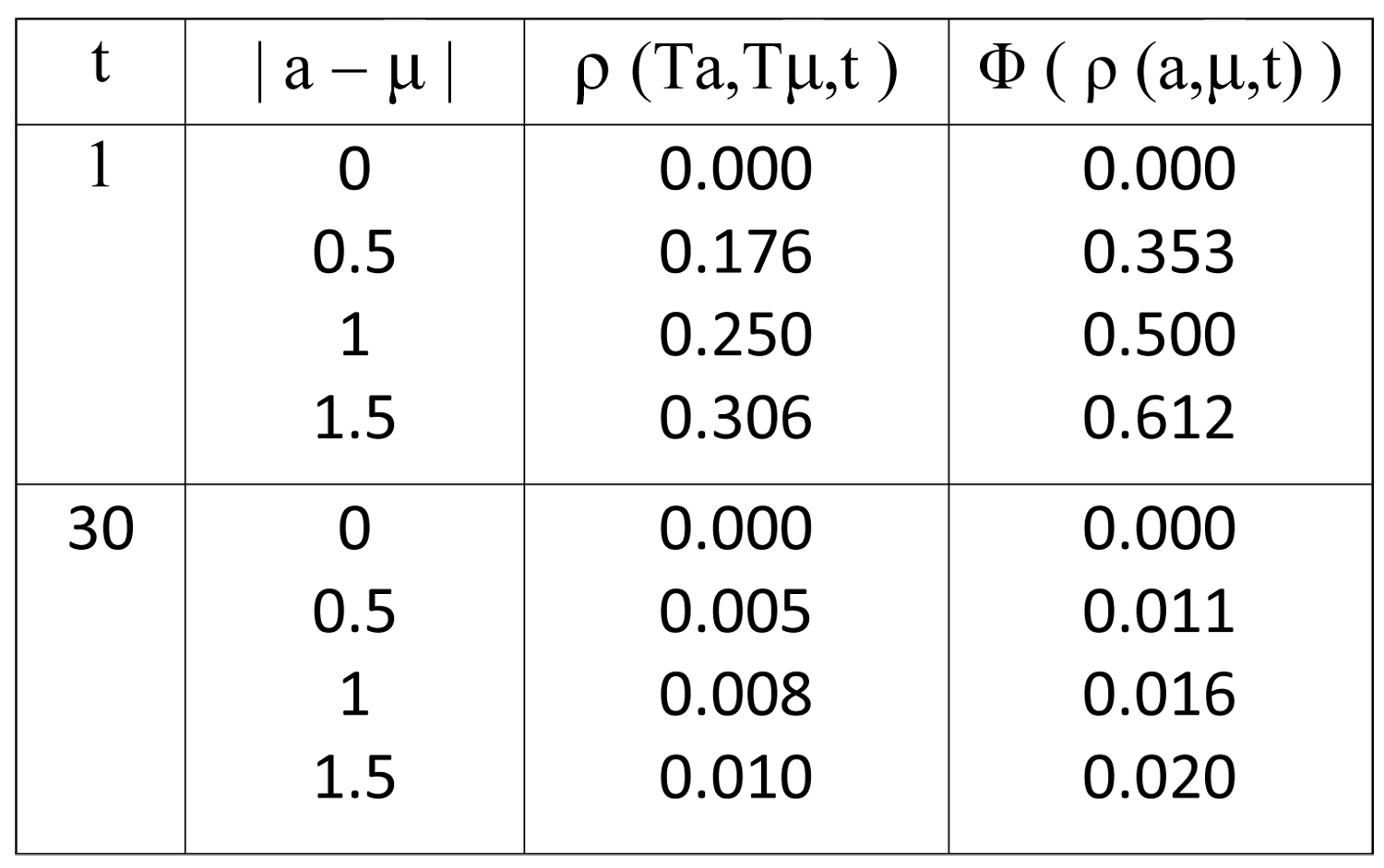}
			\caption{Variation of $ \rho ( {T} a,   T \mu, t ) $ with $ \phi ( \rho ( a, \mu, t ) ) $ as a function of $ | a - \mu | $ with fixed values of $ t = 1 $ and $  30 $}
		\end{table}
		\begin{table}[h!]
			\centering
			\includegraphics[width=10cm,height=5cm]{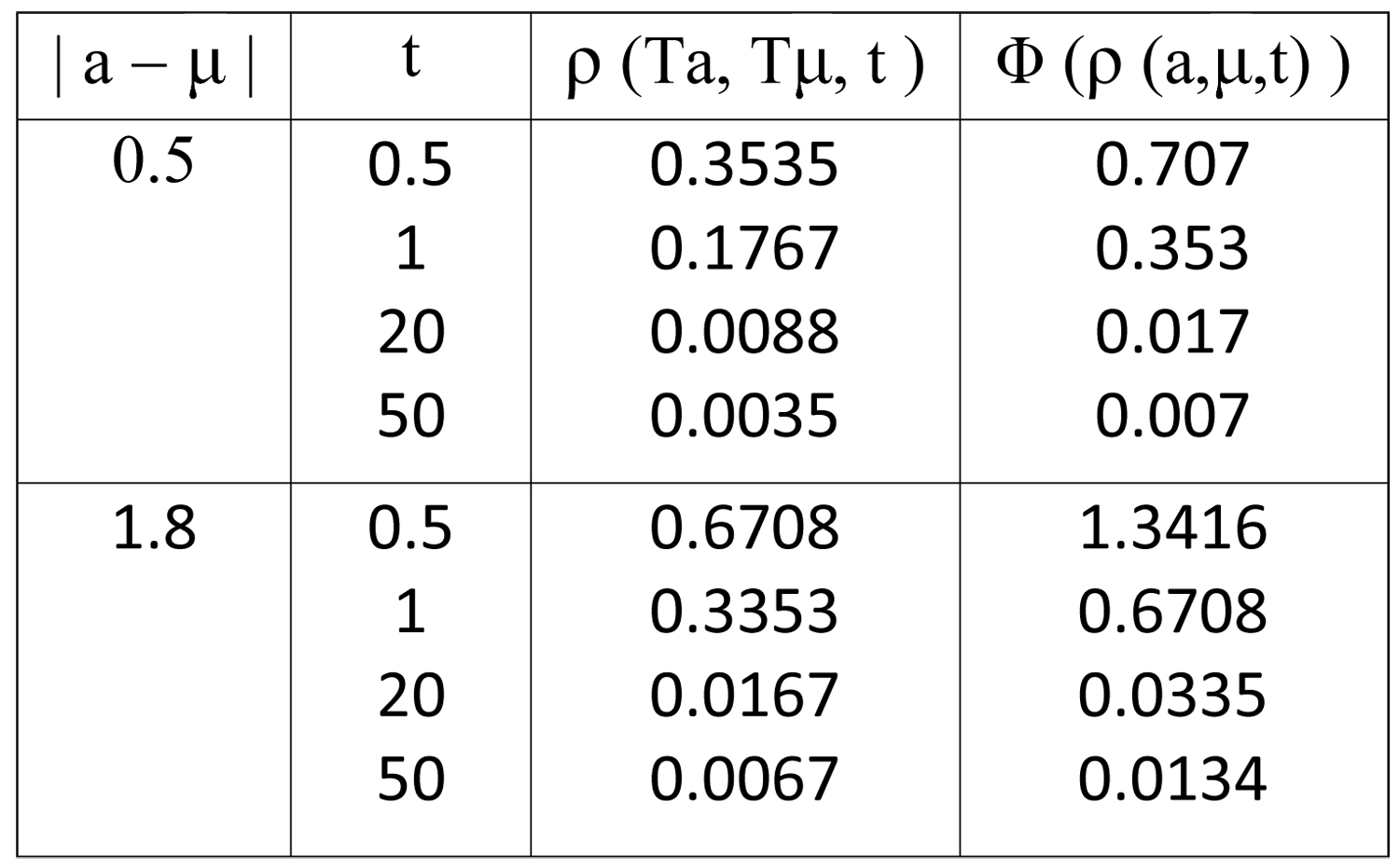}
			\caption{Variation of $ \rho ( {T} a, {T} \mu, t ) $ with $ \phi ( \rho ( a, \mu, t ) ) $ as a function of $  t  $ with fixed values of $ | a - \mu |  = 0.5 $ and $ 1.8 $}
		\end{table}
	\end{eg}

	\subsection{Application of Boyd-Wong type contraction} 
	
	Here we use the Boyd-Wong type fixed point theorem of generalized parametric metric space to  find a solution of a particular type second order initial value problem. \\ 
	For, we consider the  2nd-order IVP (\ref{B-W 1}). \\
	The differential equation (\ref{B-W 1}) is equivalent to the integral equation 
	$$ \mathcal{Y} ( \eta ) =  l_1 \cos(w \eta) + l_2 \sin (w \eta) + \int_{0}^{\eta} G ( \eta, u ) ~ \mathcal{G} ( u, \mathcal{Y} (u) )~ du, ~~~ \eta \in [ 0, S ] $$
	where $ G ( \eta, u ) $ is the Green's function defined by   
	$$ G ( \eta, u ) = \frac{1}{w} \sin ( w ( \eta - u ) ) H ( \eta - u ) $$ 
	and $  H $ is the Heaviside unit function. \\
	Next consider a complete generalized parametric metric space $ ( \Im, \rho, o ) $ where $ \Im = C [ 0, S ] $, the generalized parametric metric $ \rho $ is defined by 
	$$ \rho ( \mathcal{F}, \mathcal{G}, s ) = \dfrac{ \underset{\eta \in [ 0, S ] }{\max } | \mathcal{F} ( \eta ) - \mathcal{G} ( \eta ) | }{s} ~~~ \text{for all} ~ \mathcal{F}, \mathcal{G} 	\in \Im  ~  \text{and}  ~ s > 0 $$
	and the binary operation  `$ o $' is considered as `$ \max $'. \\

	The next theorem is an application of Boyd-Wong type fixed point theorem in generalized parametric metric space. 
	
	\begin{thm}
		Consider the differential equation (\ref{B-W 1}) where the function $ \mathcal{G} $ satisfies the following conditions:
		\begin{enumerate}[(i)]
			\item $ \mathcal{G} $ is a continuous function;
			\item there exists a non-decreasing function $ \phi : [0, \infty ) \rightarrow [0, \infty )  $ which is upper semi-continuous from right, $ \phi ( \zeta ) < \zeta ~~ \text{for all} ~~ \zeta > 0 $ and 
			\begin{equation}\label{B-W 2}
				\dfrac{1}{\zeta} | \mathcal{G} ( y, r ) - \mathcal{G} ( y, s ) | \leq w^2 \phi ( \dfrac{| r - s | }{ \zeta }) ~ ~ \text{for all}  ~ \zeta > 0, ~ y \in [ 0, S ] ~ \text{and} ~ r, s \in \mathbb{R}.
			\end{equation} 
		\end{enumerate} 
		Then there exists a unique solution for the  differential equation (\ref{B-W 1}).
	\end{thm}
	\begin{proof}
		We consider the generalized parametric metric space $ ( \Im, \rho, \max ) $ defined as above. \\
		For each $  \mathcal{Y} \in \Im $, let us define a mapping $ \mathcal{T} : \Im \rightarrow \Im $ by 
		\begin{equation}\label{B-W 3}
			(  \mathcal{T} \mathcal{Y} ) ( \mu ) = l_1 \cos ( w \mu ) + l_2 \sin ( w \mu ) + \int_{0}^{\mu} G ( \mu, u ) ~ g (u, \mathcal{Y} ( u ) ) ~ du, ~~ \text{for all} ~ \mu \in [ 0, S ].
		\end{equation}
		It very obvious that the existence of fixed points for the mapping $ \mathcal{T} $ implies the existence of solution to the equation (\ref{B-W 1}). \\
		Now for each $ \mathcal{F}_1, ~ \mathcal{F}_2 \in \Im $, $ \mathcal{Y} \in [ 0, S ] $ and $  s > 0 $, by (\ref{B-W 3}) we have
		\begin{align*}
			&	\frac{1}{s} | ( \mathcal{T} \mathcal{F}_1 ) ( \mathcal{Y} ) - ( \mathcal{T} \mathcal{F}_2 ) ( \mathcal{Y} ) |   \\
			= & 	\frac{1}{s} | \int_{0}^{\mathcal{Y}} G ( \mathcal{Y}, u ) ~ \mathcal{G} (u, \mathcal{F}_1 ( u ) ) ~ du - \int_{0}^{ \mathcal{Y} } G ( \mathcal{Y}, u ) ~ \mathcal{G} (u, \mathcal{F}_2 ( u ) ) ~ du  |  \\
			= & 	\frac{1}{s} | \int_{0}^{ \mathcal{Y} } G ( \mathcal{Y}, u ) [ \mathcal{G} (u, \mathcal{F}_1 ( u ) ) - \mathcal{G} (u, \mathcal{F}_2 ( u ) ) ] ~ du  |  \\ 
			\leq & \int_{0}^{ \mathcal{Y} } \dfrac{G ( \mathcal{Y}, u ) }{ s } | \mathcal{G} (u, \mathcal{F}_1 ( u ) ) - \mathcal{G} (u, \mathcal{F}_2 ( u ) ) | ~ du  | \\ 
			\leq & \int_{0}^{\mathcal{Y}}  G ( \mathcal{Y}, u ) ~ w^2 ~ \phi ( \dfrac{| \mathcal{F}_1 ( u ) - \mathcal{F}_2 ( u ) | }{s} )  ~ du  ~~~ (\text{by} ~ (\ref{B-W 2}))\\ 
			\leq & ~ w^2 ~ \phi ( \rho ( \mathcal{F}_1, \mathcal{F}_2, s ) )  ~ \underset{ \mathcal{Y} \in  [ 0, S ]}{\sup} \int_{0}^{ \mathcal{Y} } G ( \mathcal{Y}, u )  ~ du   \\
			= & ~ w^2 ~ \phi ( \rho ( \mathcal{F}_1, \mathcal{F}_2, s ) )  ~ \underset{ \mathcal{Y} \in  [ 0, S ]}{\sup} \int_{0}^{\mathcal{Y}} \frac{1}{w} \sin ( w ( \mathcal{Y} - u ) )  ~ du   \\
			= &  ~ \phi ( \rho ( \mathcal{F}_1, \mathcal{F}_2, s ) )  ~ \underset{ \mathcal{Y} \in  [ 0, S ]}{\sup}  ~ [ 1 - \cos  ( w \mathcal{Y} )  ]  \\
			\leq & ~ \phi (\rho ( \mathcal{F}_1, \mathcal{F}_2, s ) ). 
		\end{align*}
		This yields,
		\begin{align*}
			& ~ \underset{ \mathcal{Y} \in  [ 0, S ]}{\sup} ~	\frac{1}{s} | ( \mathcal{T} \mathcal{F}_1 ) ( \mathcal{Y} ) - ( \mathcal{T} \mathcal{F}_2 ) ( \mathcal{Y} ) | 	\leq  ~ \phi ( \rho ( \mathcal{F}_1, \mathcal{F}_2, s ) ) ~~   \text{for all} ~ s > 0 \\
			i.e  & ~~ ~ \rho (  \mathcal{T} \mathcal{F}_1,  \mathcal{T} \mathcal{F}_2, s ) 	\leq  ~ \phi ( \rho ( \mathcal{F}_1,  \mathcal{F}_2, s ) ) ~~ \text{for all} ~ s > 0.	
		\end{align*}
		Therefore by the Theorem \ref{main theorem}, $ \mathcal{T} $ is a mapping with unique fixed point, say $ \mathcal{F} $ in $ C [ 0, S ] $ and hence $ \mathcal{F} $ is the  unique solution of the system (\ref{B-W 1}) in $ C [ 0, S ] $. 
	\end{proof}

	\section{Fixed point theorems in a generalized parametric metric space over a partially ordered set }
	In this section we prove Banach type fixed point theorem in generalized parametric metric space over partially order set and then apply the result to prove the existence of solution for the first order periodic differential equation  (\ref{periodic ode1}).  
	\begin{thm} \label{sec 4 thm1}
		Let $ ( \Im, \preceq ) $ be  a poset such that each pair of $ x, \mu \in \Im $ has a lower bound and an upper bound. Assume that there exists a generalized parametric metric  $ \rho $ on $ \Im $  which makes $ ( \Im, \rho, o) $  complete  and  $`o $' be continuous. Let $ \mathcal{F}: \Im \rightarrow \Im $ be a   monotone(order-preserving or order-reversing) and continuous function satisfying 
		\begin{equation}\label{sec 4 thm1 equation 1}
			\rho ( \mathcal{F} x, \mathcal{F} \mu, s ) \leq {\kappa} ~ \rho  (x, \mu, s ) ~ ~ \text{for all} ~ ~ \mu \preceq x  ~~ \text{and} ~ ~ s > 0
		\end{equation}
		where $ 0 \leq {\kappa} < 1 $. If there exist $ \alpha \in \Im $ satisfying either $ \alpha \preceq \mathcal{F} ( \alpha ) $ or $ \alpha \succeq \mathcal{F} ( \alpha ) $ then $ \mathcal{F} $ has a unique fixed point in $ \Im $. 
	\end{thm}
	\begin{proof}
		If $ \mathcal{F} ( \alpha ) = \alpha $ then the proof is done. \\
		Suppose $ \mathcal{F} ( \alpha ) \neq \alpha $. \\
		Since $ \mathcal{F} $ is monotone on $ \Im $ and either $ \alpha \preceq \mathcal{F} ( \alpha ) $ or $ \alpha \succeq \mathcal{F} ( \alpha ) $, so we obtain a sequence in $ \Im $ satisfying  
		$$ \alpha \preceq \mathcal{F} ( \alpha ) \preceq \mathcal{F}^2 ( \alpha ) \preceq \cdots \preceq  \mathcal{F}^n ( \alpha ) \preceq \cdots ~~~ or ~~~ \alpha \succeq  \mathcal{F} ( \alpha ) \succeq  \mathcal{F}^2 ( \alpha ) \succeq  \cdots \succeq   \mathcal{F}^n ( \alpha ) \succeq  \cdots. $$
		Therefore, either  $ \mathcal{F}^{q-1} ( \alpha ) \preceq \mathcal{F}^q ( \alpha ) $ or $ \mathcal{F}^{q-1} ( \alpha ) \succeq  \mathcal{F}^q ( \alpha ), $ for any $ q \in \mathbb{N} $. \\
		Hence for all $ q \in \mathbb{N} $, from (\ref{sec 4 thm1 equation 1}) we can write 
		$$ \rho ( \mathcal{F}^{q+1} ( \alpha ), \mathcal{F}^q ( \alpha ), s ) \leq {\kappa} ~ \rho ( \mathcal{F}^{q} ( \alpha ), \mathcal{F}^{q-1} ( \alpha ), s ) ~ ~ \text{for all} ~ s > 0 $$
		where $ 0 \leq \large{\kappa} < 1 $. \\ 
		Continuing in this way we get, 
		$$ \rho ( \mathcal{F}^{q+1} ( \alpha ), \mathcal{F}^q ( \alpha ), s ) \leq {\kappa}^q ~ \rho ( \mathcal{F} ( \alpha ),  \alpha, s ) ~ ~ \text{for all} ~ s > 0 ~ \text{and}  ~ q \in \mathbb{N}. $$
		Let $ m = n + p, ~ p= 1,2,\cdots $ and $ n \in \mathbb{N} $. Then for $ m > n $ and $ s > 0 $,
		\begin{align*}
			& \rho ( \mathcal{F}^{m} ( \alpha ), \mathcal{F}^n ( \alpha ), s ) \\
			\leq & ~ \rho ( \mathcal{F}^{ n + p} ( \alpha ), \mathcal{F}^{n+1} ( \alpha ), \dfrac{s}{2} ) ~ o ~ \rho ( \mathcal{F}^{n+1} ( \alpha ), \mathcal{F}^n ( \alpha ), \dfrac{ s}{2} ) \\
			\leq & ~ \rho ( \mathcal{F}^{ n + p} ( \alpha ), \mathcal{F}^{n+2} ( \alpha ), \dfrac{s}{4} ) ~ o ~ \rho ( \mathcal{F}^{n+2} ( \alpha ), \mathcal{F}^{n+1} ( \alpha ), \dfrac{s}{4} ) ~ o ~ \rho ( \mathcal{F}^{n+1} ( \alpha ), \mathcal{F}^n ( \alpha ), \dfrac{ s}{2} )  \\
			& \cdots \hskip 20 pt \cdots \hskip 20 pt \cdots \\
			\leq & ~ \rho ( \mathcal{F}^{ n + p} ( \alpha ), \mathcal{F}^{n+p-1} ( \alpha ), \dfrac{s}{2^{p-1}} ) ~ o ~ \rho ( \mathcal{F}^{n+p-1} ( \alpha ), \mathcal{F}^{n+p-2} ( \alpha ), \dfrac{ s}{2^{p-1}} ) ~ o \cdots ~ o ~ \rho ( \mathcal{F}^{n+2} ( \alpha ), \mathcal{F}^{n+1} ( \alpha ), \dfrac{ s}{4} )   \\
			& \hskip 300 pt ~ o ~ \rho ( \mathcal{F}^{n+1} ( \alpha ), \mathcal{F}^n ( \alpha ), \dfrac{ s}{2} )  \\
			\leq & ~ \kappa^{n+p-1} \rho ( \mathcal{F} ( \alpha ),  \alpha, \dfrac{ s}{2^{p-1}} ) ~ o ~  \kappa^{n+p-2} \rho ( \mathcal{F} ( \alpha ),  \alpha, \dfrac{ s}{2^{p-1}} ) ~ o \cdots ~ o ~  \kappa^{n + 1} \rho ( \mathcal{F} ( \alpha ),  \alpha, \dfrac{ s}{4} ) ~ o ~  \kappa^{n} \rho ( \mathcal{F} ( \alpha ),  \alpha, \dfrac{ s}{2} ). 
		\end{align*}	 
		Letting $ n \rightarrow + \infty $ on the both side of the above inequality, we get 
		\begin{align*}
			& \underset{n \rightarrow + \infty}{lim} ~ \rho ( \mathcal{F}^{n+p} ( \alpha ), \mathcal{F}^n ( \alpha ), s ) \leq  \underset{n \rightarrow + \infty}{lim}  ~ \kappa^{n+p-1} \rho ( \mathcal{F} ( \alpha ),  \alpha, \dfrac{ s}{2^{p-1}} ) ~ o ~ \cdots ~ o \underset{n \rightarrow + \infty}{lim}  ~ \kappa^{n + 1} \rho ( \mathcal{F}( \alpha ),  \alpha, \dfrac{ s}{4} ) ~ o \\
			&  \hskip 262 pt  \underset{n \rightarrow + \infty}{lim}  ~ \kappa^{n} \rho ( \mathcal{F} ( \alpha ),  \alpha, \dfrac{ s}{2} ) ~~~  \text{for all} ~  s > 0 \\
			\implies & ~ \underset{n \rightarrow + \infty}{lim} ~ \rho ( \mathcal{F}^{n+p} ( \alpha ), \mathcal{F}^n ( \alpha ), s ) \leq 0 ~ o ~ 0 ~ o ~ \cdots ~ o ~ 0 ~~~  \text{for all} ~  s > 0 ~~~ (Since ~ \kappa < 1 ) \\
			\implies & ~ \underset{n \rightarrow + \infty}{lim} ~ \rho ( \mathcal{F}^{n+p} ( \alpha ), \mathcal{F}^n ( \alpha ), s )  = 0 ~~  \text{for all} ~ s > 0
		\end{align*}
		which ensures that $ \{ \mathcal{F}^n \alpha \} $ is a Cauchy sequence.  \\
		By the completeness of  $ (\Im, \rho, o ) $, there exists a converging point of  $ \{ \mathcal{F}^n \alpha \} $, say  	$ \beta  $. \\
		Next, we prove that $ \beta $ is a fixed point of $ \mathcal{F} $. \\
		Using the continuity of $ \mathcal{F} $ at $ \beta $, we have
		$$ \underset{n \rightarrow + \infty}{lim} ~ \rho ( \mathcal{F} ( \mathcal{F}^{n} ( \alpha ) ), \mathcal{F} ( \beta ),  s ) = 0 ~ ~  \text{for all} ~ s > 0 ~ ~ \text{whenever} ~ ~ \underset{n \rightarrow + \infty}{lim} ~ \rho ( \mathcal{F}^{n} ( \alpha ),  \beta, s ) = 0 ~ ~  \text{for all} ~ s > 0.  $$ 
		Now,
		\begin{align*}
			& 	\rho ( \mathcal{F} ( \beta ), \beta, s ) \leq 	\rho ( \mathcal{F} ( \beta ), \mathcal{F}^{n} (\alpha), \frac{s}{2} ) ~ o ~ \rho (  \mathcal{F}^{n} (\alpha), \beta, \frac{s}{2} ) ~ ~  \text{for all} ~ s > 0 \\
			\implies & 	\rho ( \mathcal{F} ( \beta ), \beta, s ) \leq  \underset{n \rightarrow + \infty}{lim}  	\rho ( \mathcal{F} ( \beta ), \mathcal{F}^{n} (\alpha), \frac{s}{2} )  ~ o ~ \underset{n \rightarrow + \infty}{lim} \rho (  \mathcal{F}^{n} (\alpha), \beta, \frac{s}{2} ) ~ ~  \text{for all} ~ s > 0 \\
			\implies & 	\rho ( \mathcal{F} ( \beta ), \beta, s )  = 0 ~ ~  \text{for all} ~ s > 0.
		\end{align*}		   
		This implies that $ \mathcal{F} ( \beta ) = \beta	$. \\
		It only remains to prove the uniqueness of the fixed point for $ \mathcal{F} $. \\
		We prove this by only proving that for every $ \gamma \in \Im $, $ \underset{n \rightarrow + \infty}{lim} ~  \mathcal{F}^{n} ( \gamma ) =  \beta $.  \\
		If  $ \gamma \in \Im $ be such that $ \gamma \preceq \alpha $ or $ \gamma\succeq \alpha $, then either $ \mathcal{F}^n ( \gamma ) \preceq \mathcal{F}^n ( \alpha ) $ or $ \mathcal{F}^n ( \gamma )  \succeq \mathcal{F}^n ( \alpha ), ~ ~  \text{for all} ~ n \in \mathbb{N} $, which implies  
		$$  \rho (  \mathcal{F}^{n} (\gamma), \mathcal{F}^n ( \alpha ), s ) \leq \kappa^n \rho ( \gamma, \alpha, s ) ~ ~  \text{for all} ~ s > 0. $$
		Hence,
		$$ \underset{n \rightarrow + \infty}{lim} ~ \rho (  \mathcal{F}^{n} (\gamma), \beta, s )  = 0 ~ ~  \text{for all} ~ s > 0 ~~ ~ ( \text{since } ~ \kappa < 1 ). $$
		Thus the sequence $ \{  \mathcal{F}^{n} (\gamma) \} $ also converges to $ \beta $. \\
		Next suppose $ \delta $ be an arbitrary element of $ X $ and $ \eta $ and $  \zeta $ receptively be the lower and upper bound of $ \alpha $ and $ \delta $.  \\
		Then,
		\begin{equation}\label{sec 4 thm1 equation 2}
			\eta \preceq \alpha \preceq \zeta ~~~~ \text{and} ~~~~ \eta \preceq \delta \preceq \zeta. 	
		\end{equation} 
		From (\ref{sec 4 thm1 equation 2}) it follows that either 
		\begin{equation}\label{sec 4 thm1 equation 3}
			\mathcal{F}^n (	\eta ) \preceq \mathcal{F}^n ( \alpha ) \preceq \mathcal{F}^n ( \zeta )  ~~~~ \text{and} ~~~~ \mathcal{F}^n ( \eta ) \preceq \mathcal{F}^n (\delta )  \preceq \mathcal{F}^n  ( \zeta ) ~ 
		\end{equation}
		or
		\begin{equation}\label{sec 4 thm1 equation 3.5} 
			\mathcal{F}^n (	\eta ) \preceq \mathcal{F}^n ( \alpha ) \preceq \mathcal{F}^n ( \zeta )  ~~~~ \text{and} ~~~~ \mathcal{F}^n ( \eta ) \preceq \mathcal{F}^n (\delta )  \preceq \mathcal{F}^n  ( \zeta ). 
		\end{equation}
		Again, 
		\begin{equation}\label{sec 4 thm1 equation 4}
			\underset{n \rightarrow + \infty}{lim} ~ 	\mathcal{F}^n (	\eta )  = \beta  = ~  \underset{n \rightarrow + \infty}{lim} ~ \mathcal{F}^n (	\zeta ). 
		\end{equation}
		Finally (\ref{sec 4 thm1 equation 3}), (\ref{sec 4 thm1 equation 3.5}) and (\ref{sec 4 thm1 equation 4}) together gives 
		$$   \underset{n \rightarrow + \infty}{lim} ~ 	\mathcal{F}^n (	\delta )  = \beta. $$ 
		This fulfill the requirement.
	\end{proof}

	\begin{rem}
		The Theorem \ref{sec 4 thm1} can be proved without considering the continuity of the self mapping $ f $ if we assume an additional condition.
	\end{rem}

	\begin{thm}\label{sec 4 thm2}
		Let $ ( \Im, \preceq ) $ be  a poset such that each pair of $ x, \mu \in \Im $ has a lower bound and upper bound. Suppose there exist a generalized parametric metric  $ \rho $ on $ \Im $  such that $ ( \Im, \rho, o) $ is a complete  generalized parametric metric space with continuous binary operation $`o $'. \\
		Suppose that for a  sequence $ \{ \mu_n \} $ in $ \Im $, $ \mu_n \rightarrow \mu $ as $ n \rightarrow \infty $ implies either $ \mu_n \preceq \mu $ or $ \mu_n \succeq \mu ~  \text{for all} ~  n $. Let   $ \mathcal{F} : \Im \rightarrow \Im $ be a  monotone(order-preserving or order-reversing) function satisfying  
		$$ \rho ( \mathcal{F} x, \mathcal{F} \mu, s ) \leq \kappa \rho  (x, \mu, s ) ~ ~  \text{for all} ~ ~ \mu \preceq x  ~~ \text{and} ~~ s > 0 $$ 
		where $ 0 \leq \kappa < 1 $.
		If there exist $ \alpha \in \Im $ satisfying either $ \alpha \preceq \mathcal{F} ( \alpha ) $ or $ \alpha \succeq \mathcal{F} ( \alpha )  $ then $ \mathcal{F} $ has a unique fixed point in $ \Im $. 	
	\end{thm}
	\begin{proof}
		From the proof of  Theorem \ref{sec 4 thm1}, we arrive at the converging point $ \beta $ of the sequence $ \{ \mathcal{F}^n ( \alpha ) \} $. \\
		We only prove that $ \beta $ is a fixed point of $ \mathcal{F} $. \\
		Now, since $ \{ \mathcal{F}^n ( \alpha ) \} $ converges to  $ \beta $, so $ \underset{n \rightarrow + \infty}{lim} \rho (  \mathcal{F}^{n} (\alpha), \beta, s ) = 0 ~ ~  \text{for all} ~ s > 0  $. \\
		Again $ \mathcal{F} $ is a monotone  function, so $ \{ \mathcal{F}^n ( \alpha ) \} $ is a monotone sequence converging to $ \beta $, which implies either $ \mathcal{F}^n( \alpha ) \preceq \beta $ or $  \mathcal{F}^n( \alpha ) \succeq \beta~~  \text{for all} ~ n $. \\
		Now the contraction condition gives, 
		\begin{align*}
			&	\rho ( \mathcal{F} ( \beta ), \beta, s ) \leq 	\rho ( \mathcal{F} ( \beta ), \mathcal{F}^{n+1} (\alpha), \frac{s}{2} ) ~ o ~ \rho (  \mathcal{F}^{n+1} (\alpha), \beta, \frac{s}{2} ) ~ ~  \text{for all} ~  s > 0 \\
			\implies & 	\rho ( \mathcal{F} ( \beta ), \beta, s ) \leq  \kappa ~ 	\rho (  \beta, \mathcal{F}^{n} (\alpha), \frac{s}{2} )  ~ o ~  \rho (  \mathcal{F}^{n+1} (\alpha), \beta, \frac{s}{2} ) ~ ~  \text{for all} ~ s > 0. 
		\end{align*}
		Since $ \kappa < 1 $, thus by letting $ n \rightarrow + \infty $, we have 
		$$ 	\rho ( \mathcal{F} ( \beta ), \beta, s )  \leq  0 ~ ~  \text{for all} ~ s > 0 ~~~ i.e ~~~ \rho ( \mathcal{F} ( \beta ), \beta, s )  = 0 ~ ~ \text{for all} ~ s > 0 . $$
		This implies that $ \mathcal{F} ( \beta ) = \beta	$ and hence $ \beta	$ is a fixed point for $ \mathcal{F} $.  \\
		By similar reason	as in Theorem \ref{sec 4 thm1}, we can establish the uniqueness of the fixed point. 
	\end{proof}

	Following  is an application of the above Theorem \ref{sec 4 thm1} for the existence of solution of first order differential equation.

	\begin{thm}
		Consider the problem (\ref{periodic ode1}) where $ \mathcal{F} $ is a continuous function and suppose there exist $ a > 0, ~ b > 0 $ with $ a > b $ such that for all $ r_1, r_2 \in \mathbb{R} $ with $ r_2 \geq r_1 $, $ \mathcal{F} $ satisfies 
		\begin{equation} \label{sec 4 2nd thm equ 1}
			0 \leq \{ \mathcal{F} ( y, r_2) + a r_2 \} - \{ \mathcal{F} ( y, r_1 ) + a r_1 \} \leq b ( r_2 - r_1 ).	
		\end{equation}
		Then the existence of a lower(or an upper) solution for (\ref{periodic ode1}) guarantees the existence of a unique solution of (\ref{periodic ode1}).  
	\end{thm}
	
	\begin{proof}
		Let us write the system (\ref{periodic ode1}) as
		\begin{equation}\label{periodic ode2}
			u'( y ) + a u ( y ) = \mathcal{F} ( y, u ( y ) ) + a u ( y ), ~ u ( 0 ) = u ( S ), ~ ~ y \in [0, S ].
		\end{equation}
		Then the following integral equation  
		$$ u ( y ) = \int_{0}^{S} G ( y, z ) ( \mathcal{F} ( z, u ( z ) ) + a u ( z ) ) ~ dz	$$ 
		where  $$ G ( y, z ) = \begin{cases}
			\dfrac{e^ {a ( S + z - y ) }}{e^ { aS} - 1 } ~~~~ 0 \leq z < y \leq S \\
			\dfrac{e^ {a (  z - y ) }}{ e^ { a S } - 1 } ~~~~~~~ 0 \leq y < z \leq S.
		\end{cases} $$
		is equivalent to the  system (\ref{periodic ode2}). \\ 
		Define $ \mathscr{F} : C [0, S ] \rightarrow C [0, S ]  $ by 
		\begin{equation} \label{periodic ode3}
			(\mathscr{F} u ) ( y ) = \int_{0}^{S} G ( y, z ) [ \mathcal{F} ( z, u ( z ) ) + a u ( z ) ] ~ dz ~  ~ \text{for all} ~ y \in [ 0, S ].
		\end{equation} 
		Clearly, if $ v \in C [0, S ] $ is a fixed point of $ \mathscr{F}  $ then $v $ is a solution of  (\ref{periodic ode1}). \\
		Let $ \Im = C [0, S ] $. Then $ ( \Im, \preceq ) $ is a partially ordered set whether $ `\preceq $' is defined  by:  
		$$ \mathcal{F}_1, ~ \mathcal{F}_2 \in \Im, ~ \mathcal{F}_1 \preceq \mathcal{F}_2 ~ ~ \text{if and only if} ~ ~ \mathcal{F}_1 ( y ) \leq \mathcal{F}_2 ( y )  ~  ~ \text{for all} ~  y \in [ 0, S ] .$$
		Also $ ( \Im, \rho, \max ) $ is a complete generalized parametric metric space where $ \rho $ is defined by 
		$$ \rho ( \mathcal{F}, \mathcal{G}, t ) = \dfrac{ {\underset{ y\in [0, S ] }{ \sup }} |  \mathcal{F} ( y ) - \mathcal{G} ( y ) | }{t} ~ ~  ~~\text{for all} ~ ~ \mathcal{F}, \mathcal{G} \in \Im. $$
		Again, if $ v \preceq u $ then by hypothesis (\ref{sec 4 2nd thm equ 1}), $ \mathcal{F} ( y, v ( y ) ) + a v ( y ) \preceq \mathcal{F} ( y, u ( y ) ) + a u ( y ) $ and since $ G ( x, y ) $  is non-negative for all $ ( x, y ) \in [ 0, S ] \times [ 0, S ] $, it follows that  $ \mathscr{F} $ is order-preserving on $ C [0, S ] $. \\
		Now for all $ t > 0 $ and $ u \preceq v $,
		\begin{align*}
			& ~  \dfrac{ {\underset{ y\in [ 0, S ] }{ \sup }} | (\mathscr{F} v) ( y ) - ( \mathscr{F} u )  ( y ) | }{t} ~  \\
			= & ~ \frac{1}{t} {\underset{ y\in [0, S ] }{ \sup }} \int_{0}^{S} G ( y, p ) | { \mathcal{F} ( p, v (p) ) + a v ( p ) } - { \mathcal{F} ( p, u ( p ) ) + a u ( p ) }  | ~ dp   \\
			\leq & ~  \frac{1}{t} {\underset{ y\in [0, S ] }{ \sup }} \int_{0}^{S} G ( y, p )  b | v ( p ) - u ( p ) | ~ dp ~~~ (by ~ (\ref{sec 4 2nd thm equ 1}))\\
			\leq & ~ b \rho (  v, u , t )  {\underset{ y\in [0, S ] }{ \sup }} \int_{0}^{ S } G ( y, p ) ~ dp  \\
			= & ~ b \rho (  v, u , t )  {\underset{ y\in [0, S ] }{ \sup }} [  \int_{0}^{y} \dfrac{e^{ a ( S + p - y ) }}{e^{ a S } - 1 } ~ dp + \int_{y}^{S}  \dfrac{e^{ a ( p - y ) }}{e^{ a S } - 1 } ~ dp ]   \\
			= &  \frac{b}{a}  \rho (  v, u , t ) 
		\end{align*}
		Henceforth,
		$$ \rho ( \mathscr{F} v, \mathscr{F} u, t ) \leq \frac{b}{a}  \rho (  v, u , t )  ~~ \text{for all} ~ ~ t > 0. $$ 
		Therefore  $ \mathscr{F} $ satisfies the Theorem (\ref{sec 4 thm1}) with $ \kappa = \frac{b}{a} < 1 $. \\
		Next let $ \alpha \in C [0, S ]  $ be a lower solution of (\ref{periodic ode1}). Then for all  $ y \in [ 0, S ] $ we have,
		\begin{align*}
			& \alpha' ( y ) + a \alpha ( y ) \leq \mathcal{F} ( y , \alpha ( y ) ) + a \alpha ( y ) ~~~~~  \\
			\implies & ( \alpha ( y ) e^{ a y } )' \leq [ \mathcal{F} ( y , \alpha ( y ) ) + a \alpha ( y ) ] e^{ a y } ~ ~~ ~  ~~~ (\text{multiplying both sides by} ~ e^{ a y })  \\
			\implies & \alpha ( y ) e^{ a y }  - \alpha ( 0 ) \leq \int_{0}^{y} [ \mathcal{F} ( p, \alpha ( p ) ) + a \alpha ( p ) ] e^{ a p } ~ dp. 
		\end{align*}
		Therefore, for all $   y \in [ 0, S ] $
		\begin{equation}\label{periodic ode 4}
			\alpha ( y ) e^{ a y } \leq \alpha ( 0 ) \leq \int_{0}^{y} [ \mathcal{F} ( p , \alpha ( p ) ) + a \alpha ( p ) ] e^{ a p } ~ dp. 
		\end{equation}
		Again,
		\begin{align*}
			&  \alpha ( 0  ) \leq \alpha ( S )   \\
			\implies & \alpha ( 0  ) e^{ a S } \leq \alpha ( S )  e^{ a S } \\
			& \hskip 38 pt	  \leq \alpha ( 0  )  + \int_{0}^{S} [ \mathcal{F} ( p , \alpha ( p ) ) + a \alpha ( p ) ] e^{ a p } ~ dp
		\end{align*}
		which implies 
		\begin{equation}\label{periodic ode 5}
			\alpha ( 0  ) \leq \int_{0}^{S } \dfrac{e^{ a p }}{e^{ a S } - 1 }[ \mathcal{F} ( p , \alpha ( p ) ) + a \alpha ( p ) ]  ~ dp.
		\end{equation}
		From (\ref{periodic ode 4}) and (\ref{periodic ode 5}), for all $ y \in [ 0, S ]  $ we have
		\begin{align*}
			& \alpha ( y ) e^{ a y } \leq \int_{0}^{S} \dfrac{e^{ a p }}{e^{ a S } - 1 }[ \mathcal{F} ( p , \alpha ( p ) ) + a \alpha ( p ) ]  ~ dp + \int_{0}^{y} [ \mathcal{F} ( p , \alpha ( p ) ) + a \alpha ( p ) ] e^{ a p } ~ dp ~~~~~~  \\
			\implies  &  \alpha ( y ) e^{ a y } \leq 
			\int_{0}^{y} \dfrac{e^{ a s }}{e^{ a S } - 1 }[ \mathcal{F} ( p , \alpha ( p ) ) + a \alpha ( p ) ]  ~ dp +
			\int_{y}^{S} \dfrac{e^{ a p }}{e^{ a S } - 1 }[ \mathcal{F} ( p , \alpha ( p ) ) + a \alpha ( p ) ]  ~ dp ~
			+ \\
			& \hskip 265 pt \int_{0}^{y} [ \mathcal{F} ( p , \alpha ( p ) ) + a \alpha ( p ) ] e^{ a p } ~ dp  \\
			\implies  &  \alpha ( y ) \leq \int_{0}^{y} \dfrac{e^{ a ( S + p - y ) }}{e^{ a S } - 1 }[ \mathcal{F} ( p , \alpha ( p ) ) + a \alpha ( p ) ]  ~ dp + 	\int_{y}^{S} \dfrac{e^{ a ( p - y ) }}{e^{ a S } - 1 }[ \mathcal{F} ( p , \alpha ( p ) ) + a \alpha ( p ) ]  ~ dp   \\
			\implies  &  \alpha ( y ) \leq \int_{0}^{S}  G ( y, p ) [ \mathcal{F} ( p , \alpha ( p ) ) + a \alpha ( p ) ] ~ dp \\
			\implies  &  \alpha ( y ) \leq ( \mathscr{F} \alpha ) ( y ) 
		\end{align*} 
		which implies  $ \alpha  \preceq  \mathscr{F} \alpha $. \\
		Similarly if $ \alpha \in C [0, S ]  $ be an upper  solution of (\ref{periodic ode1}) then $ ~ \alpha  \succeq  \mathscr{F} \alpha  $. \\
		Therefore from the Theorem \ref{sec 4 thm1}, $  \mathscr{F} $ has a unique fixed point in $ C [0, S ] $ and thus the existence of a unique solution of (\ref{periodic ode1}) is justified. 
	\end{proof}

	\section*{Conclusion:}
	In this manuscript we deal with the structure of some well-known metric fixed point result  in  generalized parametric metric space setting and also extend our hypothetical results of fixed point theory to partially ordered set. Here we show that these fixed point results can be applied for the existence and uniqueness criteria for solution of differential equations. \\
	This literature consists of the famous Boyd-Wong type contraction in  generalized parametric metric space and its application for the existence of solution of  second order initial value differential equation. 
	We also present an extension of `Banach contraction' in generalized parametric metric space over partially ordered set which which is also applicable  for discontinuous functions and apply it to find unique solution of first order periodic boundary value  differential equation. \\
	We hope our results will be helpful for researchers in this field for further study   and also think that there are huge scope of development in such spaces.

	\section*{Acknowledgment:} The author AD  is thankful to UGC, New Delhi, India for awarding  senior research fellowship [Grant No.1221/(CSIRNETJUNE2019)]. The authors AD and TB are  grateful to the Department of Mathematics, Siksha-Bhavana, Visva-Bharati. We are thankful to the Editor-in-Chief, Editors, and Reviewers of the journal Mathematica Slovaca for their valuable comments which are helped us to revise the manuscript in the present form.


\begin{thebibliography}{10}
		
		
		
		\bibitem{2} 
		FRÉCHET, M.:
		\textit{{Sur quelques points du calcul fonctionnel,}}  
		\newblock{Rend. del Circ. Matem. Palermo} 22 (1906), 1–72.
		
		
		
		
		\bibitem{3}
		CZERWIK, S.:
	     \textit{{Contraction mappings in b-metric spaces,}}
		\newblock{Acta Math. Univ. Ostrav.} 1 (1993), 5-11.
		
		
		
		
		
		\bibitem{4}
		SEDGHI, S.-SHOBE, N.-ALIOUCHE, A.:
		\textit{{A generalization of fixed point theorems in S-metric space,}}
		\newblock{Mat. Vesn.} 64(3) (2012), 258-266.
		
		
		
		
		\bibitem{5}  
		HUANG, L. G.-ZHANG, X.:
		\textit{{Cone metric spaces and fixed point theorems of contractive mappings,} } 
		\newblock{J. Math. Anal. Appl.} 332(2) (2007), 1468-1476. 
		
		
		
		
		\bibitem{6}
		DAS, A.-KUNDU, A.-BAG, T.:
  	    \textit{{A new approach to generalize metric functions,}}
		\newblock{Int. J. Nonlinear Anal.  Appl.} 14(3) (2023), 279–298
		
		
		
		
		\bibitem{7}
		HUSSAIN, N.-KHALEGHIZADEH, S.-SALIMI, P.-ABDOU, A. A. N.:
		\textit{{A new approach to fixed point results in triangular intuitionistic fuzzy metric spaces,} } 
		\newblock{Abstr. Appl. Anal.} 2014, (2014) Article ID 690139, 16 pages.
		
		
		
		
		\bibitem{12}
		TAŞ, N.-ÖZGÜR, N.:  
		\textit{{On parametric-metric Spaces and fixed-point type theorems for expansive mappings,}}
		\newblock{J. Math.} 2016, (2016) Article ID 4746732, 6 pages.
		
		
		
		\bibitem{13}
		DAS, A.-BAG, T.:
		\textit{{A study on parametric S-metric spaces,}}
		\newblock{Communications in Mathematics and Applications} 13(3), (2022) 921-933.
		
		
		
		\bibitem{14}
	    PRIYOBARTA, N.-ROHEN, Y.-RADENOVIC, S.:
		\textit{{Fixed point theorems on parametric A-metric space,}}
		\newblock{Am. J. Appl. Math. Stat.} 6(1), (2018) 1-5.
		
		
		
		\bibitem{15}
		TAS, N.-OZGUR, N.
		\textit{{Some fixed-point results on parametric $N_b$-metric spaces,}}
		\newblock{Commun.  Korean Math. Soc.}  33(3), (2018) 943-960.
		
		
		
		\bibitem{16}
		CETKIN, V.:
	\textit{{Parametric 2-metric spaces and some fixed point results}},
		\newblock{New Trend Math. Sci.} 7(4), (2019) 503-511.
		
		
		
		\bibitem{1}
		DAS, A.-BAG, T.
	\textit{{A generalization to parametric metric spaces},}
		\newblock {Int. J. Nonlinear Anal.  Appl.} 14(1), (2023) 229-244.
		
		
		
		
		\bibitem{8}
		BOYD, D. W.-WONG, J. S. W.:
	\textit{	{On nonlinear contractions}},
		\newblock{Proceedings of the American Mathematical Society} 20(2), (1969) 458-464. 
		
		
		
			\bibitem{new ref}
	    SOM, S.- BERA, A.-DEY, L.K.:
		\textit{	{ Some remarks on the metrizability of $\mathscr{F}$-metric spaces}},
		\newblock{J. Fixed Point Theory Appl.} 22(17), (2020). 
		
		
		
		
		\bibitem{9}
		RAN, A.-REURINGS, M.:  
		\textit{{A fixed point theorem in partially ordered sets and some applications to matrix equations}},
		\newblock {Proceedings of the American Mathematical Society} 132(5), (2003) 1435-1443.  
		
		
		
		
		
		
		
		
		\bibitem{10}
		BANACH, S.
		{\textit{Sur les operations dans las ensembles abstraits et leur application aux 
			equations integrales},}  
		\newblock{Fundam. Math.} 3(1), (1922) 133-181. 
		
		
		
		\bibitem{11}
		LADDE, G. S.-LAKSHMIKANTHAM, V.-VATSALA, A. S.:
	\textit{	{Monotone iterative techniques for nonlinear differential equations}},
		\newblock{Pitman Publishing}, 27 (1985). 
		
		
		
	
		
	\end{thebibliography}
\end{document}